\documentclass[12pt,reqno]{amsart}

\usepackage{amssymb}
\usepackage{amsmath}
\usepackage{latexsym}
\usepackage{amsthm}

\theoremstyle{plain}
\newtheorem{theorem}{Theorem}[section]

\newtheorem{lemma}{Lemma}[section]
\newtheorem{remark}{\it Remark}[section]
\theoremstyle{definition}

\newcommand{\M}{\mathrm{M}}

\newcommand{\rot}{\mathop\mathrm{rot}}

\renewcommand{\div}{\mathop\mathrm{div}}
\newcommand{\vphi}{\varphi}
\newcommand{\ve}{\varepsilon}

\title
[  Analyticity of the solutions of the 2D damped-driven NS system]
 {On the domain of analyticity  and small scales for the solutions of
 the  damped-driven 2D Navier--Stokes equations
 \footnote{Date: January 18, 2007.}}
\author[A.A. Ilyin and E.S. Titi]{}

\begin{document}

\maketitle

\centerline{\scshape  Alexei A. Ilyin\footnote{\noindent Keldysh
Institute of
Applied Mathematics, Russian Academy of Sciences, Miusskaya Sq. 4, 125047
Moscow, Russia, E-mail: ilyin@keldysh.ru } and
Edriss S. Titi\footnote
{Department of Mathematics and Department of Mechanical
and Aerospace Engineering, University of California,
 Irvine, California 92697, USA, E-mail: etiti@math.uci.edu. Also:
  Department of Computer Science and Applied Mathematics,
  Weizmann Institute of Science,
  P.O. Box 26,
  Rehovot, 76100, Israel,
  E-mail: edriss.titi@weizmann.ac.il
}}

\medskip

\bigskip
\begin{quote}{\normalfont\fontsize{8}{10}\selectfont
{\bfseries Abstract.} We obtain a logarithmically sharp estimate for
the space-analyticity radius of the solutions of the damped-driven
2D Navier--Stokes equations with periodic boundary conditions and
relate this to the small scales in this system. This system is
inspired by the  Stommel--Charney barotropic ocean circulation
model.
\medskip

\textbf{Key words:} Analyticity, Gevrey regularity,
 Navier--Stokes equations, dissipative length scales,
 Stommel--Charney model.
\medskip

\textbf{AMS  subject classification:} 35B41, 35Q30, 37L30.
\par}
\end{quote}

\setcounter{equation}{0}
\section{Introduction}\label{S:Intro}

It was shown in \cite{F-T-Gev} (see also~\cite{Chae}, \cite{FMRT}) that the solutions
of the 2D Navier--Stokes equations with periodic boundary conditions
belong to the Gevrey class of analytic functions
(if the forcing term does).
Using the Gevrey regularity approach the following estimate for the
spatial analyticity radius for the solutions
that lie on the global attractor (or are near it) was obtained
\begin{equation}\label{estFT}
    l_a\ge\frac{c|\Omega|^{1/2}}{G^2\log G},
\end{equation}
where $G=\|f\|_{L_2}|\Omega|/\nu^2$ is the Grashof number and
$|\Omega|=L^2/\gamma$ is the area of the periodic domain
 $\Omega=[0,L/\gamma]\times[0,L]$, $\gamma\le1$.

 Therefore, the
Fourier coefficients $\hat u_k$ are exponentially small for
$|k|\gg L/l_a$,
and $l_a$
naturally forms a lower bound for the small dissipative length
scale for the system (see, for instance,~\cite{D-T}).

There are other ways of estimating the dissipative small length scale for the
Navier--Stokes system, for instance, in terms of the dimension
of the global attractor \cite{BV}, \cite{CF88}, \cite{CFT85}, \cite{FMRT}, \cite{T}.
The  Hausdorff and fractal dimensions  of the global attractor
satisfy the following estimate \cite{CFT} (see also~\cite{CF88},
\cite{T}):
$$
\dim_F\mathcal{A}\le c_1 G^{2/3}(\log(1+G))^{1/3},\qquad
c_1=c_1(\gamma)
$$
which has been shown in \cite{Liu} (following ideas of~\cite{BV})
 to be logarithmically sharp.

If we accept the point of view that the small length scale can be
defined as follows (see \cite{CFT85}, \cite{FMRT}, \cite{Rob}, \cite{T})
\begin{equation}\label{viadim}
    l_f\sim\left(\frac{|\Omega|}{\dim_F\mathcal{A}}\right)^{1/2},
\end{equation}
then up to logarithmic correction we have
\begin{equation}\label{heurest}
    l_f\sim\frac{|\Omega|^{1/2}}{G^{1/3}}\,.
\end{equation}
This heuristic estimate for the small length scale is
probably the best one can hope for since it matches, up to logarithmic term,
the physically asserted estimates for the enstrophy dissipation length
scale~\cite{Kr} . We also observe that the
estimate~(\ref{heurest}) is
extensive, that is, independent of the size of the spatial domain
provided that its shape is fixed.

Another rigorous definition of the small length scale can be given in
terms of the number of determining modes, nodes, or volume elements
 (see  \cite{FMRT}, \cite{F-P}, \cite{F-T-Nodes},
\cite{Jones-Titi2} and the references therein).
 It was shown that if $N$ is sufficiently large and
$N$ equal squares of size $l_{dn}$ tile the periodic spatial domain,
then any collection of points (one in each square) are determining
for the long time dynamics of the 2D Navier--Stokes system. The best  to date estimate for $N$
was obtained in \cite{Jones-Titi2}:
$$
N\le c_2G,
$$
where $c_2=c_2(\gamma)$ depends only on the aspect ratio $\gamma\le1$.
(An explicit estimate for $c_2$ was obtained in \cite{I-T3}:
$c_2(\gamma)=(68/(\gamma\pi))^{1/2}$.)

Therefore the small length scale defined in terms of the lattice
of determining nodes satisfies
\begin{equation}\label{l-nodes}
l_{dn}\ge c_2(\gamma)^{-1/2}\frac{|\Omega|^{1/2}}{G^{1/2}}\,.
\end{equation}
We observe that this estimate is not extensive, that is, $l_{dn}$ scales
like $\lambda^{-1/2}$ if $\Omega$ is replaced by $\lambda\Omega$,
$\lambda>0$.

We point out here that for the 2D Navier--Stokes system with analytic
forcing the results of~\cite{F-K-R}, \cite{F-R} provide the existence
of a finite number $N$ of instantaneously determining nodes comparable
with the fractal dimension of the attractor. These nodes, however, can
be chosen arbitrarily (up to a subset of $\Omega^N$ with $2N$-dimensional Lebesgue
measure zero) and therefore do not naturally define a regular
lattice of determining nodes.

The best to date estimate for the analyticity radius of the solutions of the
Navier--Stokes equations with analytic forcing term $f$ was obtained in
\cite{Kuk}:
\begin{equation}\label{l-a-K}
l_a\ge c_3(\gamma)\frac{|\Omega|^{1/2}}{G^{1/2}(1+\log G)^{1/4}}\,.
\end{equation}
Relating the radius of analyticity to the dissipative small length scale
(see also~\cite{H-K-R} in this regard) we
note that up to a logarithmic correction the estimate~(\ref{l-a-K})
coincides with
(\ref{l-nodes}), but both are worse than~(\ref{viadim}),
where the latter coincides, as we have already pointed out,
with the physically asserted estimate of~\cite{Kr}.

In this paper we focus on the 2D space periodic
Navier--Stokes system with damping
\begin{equation}\label{DNS}
\begin{aligned}
\partial_tu+\sum_{i=1}^2u^i\partial_iu&=
-\mu u+\nu\Delta\,u-\nabla\,p +f,\\
\div u&=0.
\end{aligned}
\end{equation}
By adding the Coriolis forcing term to~(\ref{DNS}) one obtains
the well-known Stommel--Charney barotropic model of ocean
circulation~\cite{Charney}, \cite{D-F}, \cite{P}, \cite{Stommel}.
Here the damping $\mu u$ represents the
Rayleigh friction term and $f$ is the wind stress.
For an analytical study of this system see, for instance
\cite{BCT}, \cite{Hauk}, \cite{I91},  \cite{W}, and the references therein.
In a follow up work
we will be studying the effect of adding rotation (Coriolis parameter)
on the size of small scales and the complexity of the dynamics of~(\ref{DNS}).
Therefore, we will focus in this work on the system~(\ref{DNS}).
  We  also  point out that in this geophysical context the viscosity
plays a  much smaller role in the  mechanism of dissipating energy than the Rayleigh
friction. That is why in this work the friction coefficient $\mu>0$
will be fixed and we consider the system at the limit when $\nu\to0^+$.

Sharp estimates (as $\nu\to0$) for the Hausdorff and the fractal dimensions of the global attractor of the
system~(\ref{DNS}) were first obtained in the case of the square-shaped domain
 in \cite{I-M-T} ($\gamma=1$). Then the case of  an elongated domain was studied in~\cite{I-T6} ($\gamma\to0$), where it was shown that
\begin{equation}\label{dimD}
\dim_F\mathcal{A}\le c_4D,\qquad
D=\frac{\|\rot f\|_\infty|\Omega|}{\mu\nu}\,,
\end{equation}
where $c_4$ is an absolute constant ($c_4\le12$). This  estimate
 is sharp as both $\nu\to0$ and $\gamma\to0$. Therefore the small
length scale defined as in~(\ref{viadim}) is of the order of
\begin{equation}\label{lfdamped}
l_f\sim\left(\frac{|\Omega|}{\dim_F\mathcal{A}}\right)^{1/2}
\sim\,\left(\frac{\mu\nu}{\|\rot f\|_\infty}\right)^{1/2}
\sim\,\frac{|\Omega|^{1/2}}{D^{1/2}}\,.
\end{equation}
This heuristic estimate is, in fact, a rigorous bound for the small
length scale expressed in terms of the number of determining
modes and nodes~\cite{I-T3}:
\begin{equation}\label{dndamped}
l_{dn}= c_5\left(\frac{|\Omega|}{D}\right)^{1/2}=c_5
\left(\frac{\mu\nu}{\|\rot f\|_\infty
}\right)^{1/2}, \qquad c_5=68^{1/4}.
\end{equation}
This means that any lattice of points in $\Omega$ at a typical distance
$l\le l_{dn}$ is determining.

The main result of this paper is in showing that the analyticity radius
$l_a$ of the solutions of the damped-driven Navier--Stokes
system~(\ref{DNS}) lying on the global attractor is bounded from below
and satisfies the estimate:
\begin{equation}\label{la}
l_a\ge \frac{c|\Omega|^{1/2}}{D^{1/2}(1+\log D)^{1/2}}\,,
\end{equation}
which up to a logarithmic correction agrees both with the
smallest scale  estimate~(\ref{lfdamped}) and the rigorously
defined typical distance between the determining nodes~(\ref{dndamped}).

It is worth mentioning that this point of view of relating the radius of analyticity
of solutions on the Navier--Stokes equations to small scales in turbulence
was also presented in~\cite{H-K-R}.

This paper is organized as follows.
In section~\ref{S:Gevrey} we employ the Gevrey--Hilbert space technique
of \cite{F-T-Gev} to derive a lower bound for the radius of analyticity of
the order
\begin{equation}\label{Gevrey-D}
\frac{c|\Omega|^{1/2}}{D^2\log D}\,.
\end{equation}
This bound considerably improves,  for a fixed $\mu>0$, the lower bound~(\ref{estFT})
for the classical Navier--Stokes system
as $\nu\to0^+$ (see also Remark~\ref{R:with-mu}).
Let us remark that as an alternative to the Gevrey regularity technique
for estimating small scales one can apply the
ladder estimates approach presented in~\cite{D-G}
to obtain estimates for the small scales in~(\ref{DNS})
(see also~\cite{Gi-Ti}).

In section~\ref{S:Kuk} the estimate~(\ref{la}) is proved
for the system~(\ref{DNS})  following~\cite{Kuk}.

\setcounter{equation}{0}
\section{Gevrey regularity of
the damped Navier--Stokes system}\label{S:Gevrey}
As usual (see, for instance,  \cite{BV},\cite{CF88},\cite{Lad},\cite{TNS}), we write~(\ref{DNS}) as an evolution equation
in the Hilbert space $H$ which is the closed subspace of solenoidal
vectors in $(L_2(\Omega))^2$ with zero average over the torus
$\Omega=[0,L/\gamma]\times[0,L]$:
\begin{equation}\label{FDNS}
\partial_tu+B(u,u)+\nu Au+\mu u=f, \qquad u(0)=u_0.
\end{equation}
Here $A=-P\Delta$ is the Stokes operator with eigenvalues
$0<\lambda_1\le\lambda_2\le\dots$,
$B(u,v)=P\bigl(\sum_{i=1}^2u^i\partial_iv\bigr)$ is the nonlinear term,
$f=Pf\in H$, and $P:(L_2(\Omega))^2\to H$.

We restrict ourselves to the case $\gamma=1$ and, in addition, assume   that $\Omega=[0,2\pi]^2$ (this simplifies the Fourier series below).
The case of the square-shaped domain $\Omega=[0,L]^2$ reduces to this case by scaling.
Furthermore, any domain with aspect ratio $\gamma<1$ can be treated in
the similar way, the absolute dimensionless constants $c_1, c_2,\dots$ below
will then depend on $\gamma$, however.

A vector field
$u\in H$ has the Fourier series expansion
$$
u=\sum_{j\in\mathbb{Z}^2} u_je^{ij\cdot x},
\quad u_j\in\mathbb{C}^2,
\quad u_{-j}=\bar u_j,
\quad u_j\cdot j=0,
\quad u_0=0,
$$
and
$$
\|u\|^2=\|u\|^2_{L_2}=(2\pi)^2\sum_{j\in\mathbb{Z}^2} |u_j|^2.
$$
The eigenvalues of the Stokes operator $A$ are the numbers $|j|^2$,
and the domain of its powers  is the set of vector
functions $u$ such that
$$
(2\pi)^2\sum_{j\in\mathbb{Z}^2}|j|^{4\alpha} |u_j|^2=
\|A^\alpha u\|^2<\infty.
$$
For $\tau, s>0$ we define the Gevrey space $D(e^{\tau A^s})$
of functions $u$ satisfying
\begin{equation}\label{def-Gev}
(2\pi)^2\sum_{j\in\mathbb{Z}^2} e^{2\tau|j|^{2s}} |u_j|^2=
\|e^{\tau A^s}u\|^2<\infty.
\end{equation}
We suppose that the forcing term $f$ belongs to the Gevrey space
of analytic functions
\begin{equation}\label{f-Gev}
f\in D(e^{\sigma_1A^{1/2}}A^{1/2}),
\ \text{so that}\
(2\pi)^2\sum_{j\in\mathbb{Z}^2}|j|^2e^{2\sigma_1|j|} |u_j|^2=
\|e^{\sigma_1A^{1/2}}A^{1/2}f\|<\infty
\end{equation}
for some $\sigma_1>0$. We set
$$
\vphi(t)=\min(\nu\lambda_1^{1/2}t,\sigma_1).
$$
The norm and the scalar product in $D(e^{\vphi(t)A^{1/2}})$
are denoted by $\|\cdot\|_\vphi$ and $(\cdot,\cdot)_\vphi$,
respectively.

We assume that $u_0\in D(A^{1/2})$ and take the scalar product
of~(\ref{FDNS}) and $Au$ in $D(e^{\vphi(t)A^{1/2}})$ for
sufficiently small $t\le\sigma_1/(\nu\lambda_1^{1/2})$.
Since
$$
\bigl(e^{\vphi(t)A^{1/2}}\partial_tu(t),e^{\vphi(t)A^{1/2}}u(t)\bigr)=
\frac12\partial_t\|A^{1/2}u(t)\|^2_\vphi-
\nu\lambda_1^{1/2}(Au(t),A^{1/2}u(t))_\vphi,
$$
we obtain
\begin{equation}\label{scal-pr}
\frac12\partial_t\|A^{1/2}u\|_\vphi^2+\nu\|Au\|^2_\vphi+\mu\|A^{1/2}u\|^2_\vphi=
(B(u,u)Au)_\vphi+\nu\lambda_1^{1/2}(Au,A^{1/2})_\vphi+
(A^{1/2}f,A^{1/2}u)_\vphi.
\end{equation}
Next we use the key estimate (see \cite{F-T-Gev}, \cite{FMRT}, \cite{Titi})
for the nonlinear term in Gevrey spaces
$$
(B(u,u),Au)_\vphi\le c_1\|A^{1/2}u\|^2_\vphi\|Au\|_\vphi
\biggl(1+
\log\frac{\|Au\|^2_\vphi}{\lambda_1\|A^{1/2}u\|^2_\vphi}\biggr)^{1/2}
$$
and use Young's inequality for this estimate and for the last two terms
in~(\ref{scal-pr}):
$$
\aligned
\partial_t\|A^{1/2}u\|_\vphi^2+&\nu\|Au\|^2_\vphi\le\\&\le
\frac{2c_1^2}\nu\|A^{1/2}u\|^4_\vphi
\biggl(1+\log\frac{\|Au\|^2_\vphi}{\lambda_1\|A^{1/2}u\|^2_\vphi}\biggr)
+2\nu\lambda_1\|A^{1/2}u\|^2_\vphi+\frac{\|A^{1/2}f\|^2_\vphi}{2\mu}\le\\
&\le
\frac{c_2}\nu\|A^{1/2}u\|^4_\vphi
\biggl(1+\log\frac{\|Au\|^2_\vphi}{\lambda_1\|A^{1/2}u\|^2_\vphi}\biggr)
+\nu^3\lambda_1^2+\frac{\|A^{1/2}f\|^2_\vphi}{2\mu},
\endaligned
$$
where $c_2=2c_1^2+1$. Next, using the inequality
$
-\alpha z+\beta(1+\log z)\le\beta\log\beta/\alpha
$
(see  \cite{FMRT}, \cite{FMT}), we find
$$
\aligned
-\nu\|Au\|_\vphi^2+
\frac{c_2}\nu\|A^{1/2}u\|^4_\vphi
\biggl(1+\log\frac{\|Au\|^2_\vphi}{\lambda_1\|A^{1/2}u\|^2_\vphi}\biggr)\le
\frac{c_2}\nu\|A^{1/2}u\|^4_\vphi
\log\frac{c_2\|A^{1/2}u\|^2_\vphi}{\lambda_1\nu^2},
\endaligned
$$
and obtain the differential inequality
$$
\partial_t\|A^{1/2}u\|_\vphi^2\le
\frac{c_2}\nu\|A^{1/2}u\|^4_\vphi
\log\frac{c_2\|A^{1/2}u\|^2_\vphi}{\lambda_1\nu^2}
+\nu^3\lambda_1^2+\frac{\|A^{1/2}f\|^2_\vphi}{2\mu}\,.
$$
Hence the function
$$
y(t)=\frac{c_2\|A^{1/2}u\|_\vphi^2}{\lambda_1\nu^2}+
\frac{\|A^{1/2}f\|_{\sigma_1}^2}{\lambda_1\nu^{3/2}\mu^{1/2}}
+e,
$$
where $\ln e=1$, satisfies
$$
\partial_ty(t)\le\nu\lambda_1c_3y^2\log y, \qquad c_3=\max(1,c_2/2).
$$
Therefore
$y(t)\le2y(0)$ for as long as
$$t\le(2\nu\lambda_1c_3y(0)\log2y(0))^{-1}.$$
In other words,
$$
\|A^{1/2}u\|_\vphi^2\le2\|A^{1/2}u_0\|^2+
c_4(\nu/\mu)^{1/2}\|A^{1/2}f\|_{\sigma_1}+c_4\lambda_1\nu^2,
\qquad c_4=e/c_2,
$$
as long as $0\le t\le T^*(\|A^{1/2}u_0\|)$, where
$$
\aligned
T^*(\|A^{1/2}u_0\|)=\frac1{2c_3\nu\lambda_1
\left(\frac{c_2\|A^{1/2}u_0\|^2}{\lambda_1\nu^2}
+\frac{\|A^{1/2}f\|_{\sigma_1}}{\lambda_1\nu^{3/2}\mu^{1/2}}+e
\right)\log\left(2\left(\frac{c_2\|A^{1/2}u_0\|^2}{\lambda_1\nu^2}
+\frac{\|A^{1/2}f\|_{\sigma_1}}{\lambda_1\nu^{3/2}\mu^{1/2}}+e
\right)\right)}.
\endaligned
$$

We now observe (see Lemma~\ref{L:Est-for-rot}) that on the global
attractor or in the absorbing ball we have, respectively,
$$
\|A^{1/2}u(t)\|\le\frac{\|A^{1/2}f\|}\mu,\quad t\in\mathbb{R},
\quad
\|A^{1/2}u(t)\|\le2\frac{\|A^{1/2}f\|}\mu,\quad t\ge T_0(\|A^{1/2}u_0\|).
$$
Therefore we have the following lower bound for $T^*$:
$$
T^*\ge
c_5\left[\nu\lambda_1
\left(\frac{\|A^{1/2}f\|^2}{\lambda_1\nu^2\mu^2}
+\frac{\|A^{1/2}f\|_{\sigma_1}}{\lambda_1\nu^{3/2}\mu^{1/2}}+1
\right)\log\left(\frac{\|A^{1/2}f\|^2}{\lambda_1\nu^2\mu^2}
+\frac{\|A^{1/2}f\|_{\sigma_1}}{\lambda_1\nu^{3/2}\mu^{1/2}}+1
\right)
\right]^{-1}
$$
In the limit $\nu\to0^+$ we have
$$
\frac{\|A^{1/2}f\|^2}{\lambda_1\nu^2\mu^2}\gg
\frac{\|A^{1/2}f\|_{\sigma_1}}{\lambda_1\nu^{3/2}\mu^{1/2}},
$$
and we can write the lower bound for $T^*$ as follows
$$
T^*\ge c_6\left[\nu\lambda_1D^2\log D\right]^{-1},
$$
where
$$
\frac{\|A^{1/2}f\|}{\lambda_1^{1/2}\nu\mu}=
\frac{\|\rot f\||\Omega|^{1/2}}{2\pi\nu\mu}\le\frac1{2\pi}D,
\quad \text{where}\quad
D=\frac{\|\rot f\|_\infty|\Omega|}{\nu\mu}\,.
$$
In terms of the analyticity radius $l_a$ the lower bound for
$T^*$ takes the form
$$
l_a\ge\frac{c_7|\Omega|^{1/2}}{D^2\log D}\,.
$$
Thus, we have proved the following theorem.
\begin{theorem}\label{T:Gevrey}
Suppose that $f\in D(A^{1/2}e^{\sigma_1A^{1/2}})$ for some
$\sigma_1>0$. Then a solution $u$ lying on the global attractor
$\mathcal{A}$ is analytic with analyticity radius
$$
l_a\ge\min\left(\frac{c_7|\Omega|^{1/2}}{(D^2+D_1+1)\log(D^2+D_1+1)}\ ,
 \ \sigma_1\right),
$$
where
$$
D=\frac{\|\rot f\|_\infty|\Omega|}{\nu\mu}\,,\qquad
D_1=\frac{\|A^{1/2}f\|_{\sigma_1}}{\lambda_1\nu^{3/2}\mu^{1/2}}.
$$
Moreover,
\begin{equation}\label{Gev-with-mu}
l_a\ge\frac{c_8|\Omega|^{1/2}}{D^2\log D}\,\quad
\text{as}\quad\nu\to0^+.
\end{equation}
The constants $c_7$ and $c_8$ depend only on the aspect ratio
of the periodic domain $\Omega$.
\end{theorem}

\begin{remark}\label{R:with-mu}
\rm{
We observe that the estimate~(\ref{Gev-with-mu})
for the system~(\ref{DNS}) is of the order $\nu^{-2}\log(1/\nu)$
as far as the dependence on $\nu\to0^+$ is concerned,
while the estimate~(\ref{estFT}) for the classical Navier--Stokes
system is, in this respect much larger; namely, is of the
order~$\nu^{-4}\log(1/\nu)$.

However, the estimate~(\ref{Gev-with-mu}) is not sharp and will
be improved in the next section. As has been demonstrated in~\cite{O-T}
the Gevrey--Hilbert space technique
 does not always provide sharp estimates for the radius of analyticity.
 The mechanism explaining this  has been reported in \cite{O-T} by means of an explicitly solvable
 model equation.

}
\end{remark}

\setcounter{equation}{0}
\section{Sharper bounds}\label{S:Kuk}

In this section  we obtain sharper lower bounds for the analyticity radius
$l_a$. This is achieved by combining the $\nu$-independent estimate for the
vorticity contained in the following lemma and the $L_p$-technique
developed in~\cite{Gr-K}, \cite{Kuk} for the uniform analyticity radius of
the solutions of the Navier--Stokes equations. We observe that similar technique
has  been earlier established in~\cite{B-B} for studying the analyticity
of the Euler equations.

Applying the operator $\rot$ to~(\ref{DNS}) we obtain the well-known
scalar vorticity equation
\begin{equation}\label{VE}
\partial_t\omega+u\cdot\nabla\omega=\nu\Delta\omega-\mu \omega+F,
\end{equation}
where $\omega=\rot u$, $F=\rot f$, $u=\nabla^\perp\Delta^{-1}\omega$,
 so that
$u\cdot\nabla\omega=\nabla^\perp\Delta^{-1}\omega\cdot\nabla\omega=
J(\Delta^{-1}\omega,\omega)$,
and $\nabla^\perp\,=(-\partial_2\ ,\partial_1\ )$,
$J(a,b)=\nabla^\perp a\cdot \nabla b$.

\begin{lemma}\label{L:Est-for-rot} {\rm(See~\cite{I-T3}.)}
The solutions $u(t)$ lying on the global attractor $\mathcal{A}$
satisfy the following bound:
\begin{equation}\label{rot-infty}
  \|\omega(t)\|_{L_{2k}}\le
    \frac{\|\rot f\|_{L_{2k}}}{\mu}\ ,
\qquad t\in\mathbb{R},
\end{equation}
where $1\le k\le\infty$.
\end{lemma}
\begin{proof}
We use the vorticity equation~(\ref{VE})
and take the scalar product with $\omega^{2k-1}$, where $k\ge1$ is
integer,
and use the identity $$(J(\psi,\vphi),\vphi^{2k-1})=
(2k)^{-1}\int J(\psi,\vphi^{2k})dx=
(2k)^{-1}\int\div(\vphi^{2k}\nabla^\perp\psi)dx=0.$$
We obtain
$$
\aligned
\|\omega\|_{L_{2k}}^{2k-1}\partial_t\|\omega\|_{L_{2k}}+
(2k-1)\nu\int|\nabla\,\omega|^2\omega^{2k-2}dx+
\mu\|\omega\|_{L_{2k}}^{2k}=\\=(\rot f,\omega^{2k-1})\le\
\|\rot f\|_{L_{2k}}\|\omega\|_{L_{2k}}^{2k-1}.
\endaligned
$$
Hence, by Gronwall's inequality
$$
\|\omega(t)\|_{L_{2k}}\le\|\omega(0)\|_{L_{2k}}e^{-\mu t}+
\mu^{-1}\|\rot f\|_{L_{2k}}(1-e^{-\mu t}),
$$
and passing to the limit as $k\to\infty$ we find
$$
\|\omega(t)\|_\infty\le\|\omega(0)\|_\infty e^{-\mu t}+
\mu^{-1}\|\rot f\|_\infty(1-e^{-\mu t}).
$$
Now, we let $t\to\infty$ in the above inequalities and  obtain
$$
\limsup_{t\to\infty}\|\omega(t)\|_{L_{2k}}\le
\frac{\|\rot f\|_{L_{2k}}}\mu\ ,
\quad
1\le k\le\infty,
$$
which gives (\ref{rot-infty}) since the solutions lying
on the attractor are bounded for $t\in\mathbb{R}$.
\end{proof}

As before we consider the square-shaped domain $\Omega=[0,L]^2$ and
it is now convenient to write~(\ref{DNS})
in dimensionless form.
We introduce dimensionless variables $x'$, $t'$, $u'$ and $p'$
by setting
$$
x=Lx',\quad t=(L^2/\nu)t',\quad u=(\nu/L)u',\quad p=(\nu^2/L^2)p',
\quad \mu=(\nu/L^2)\mu'.
$$
We obtain
\begin{equation}\label{DDNS}
\begin{aligned}
\partial_{t'}u'+\sum_{i=1}^2{u'}^i\partial'_iu'&=
-\mu' u+\Delta'\,u'-\nabla'\,p' +f',\\
{\div}' u'&=0,
\end{aligned}
\end{equation}
where $x'\in\Omega'=[0,1]^2$, $f'=(L^3/\nu^2)f$. Accordingly,
the dimensionless form of~(\ref{VE}) is as follows (we omit the primes):
\begin{equation}\label{DVE}
\partial_t\omega+u\cdot\nabla\omega=\Delta\omega-\mu \omega+F.
\end{equation}

\begin{remark}\label{Rm:Dimless-est}
{\rm
For dimensionless variables $u'$ and $\omega'$ the
estimate~(\ref{rot-infty}) with $k=\infty$ takes the form
\begin{equation}\label{Dimless-est}
\|\omega'\|_\infty=\|{\rot}'u'\|_\infty\le
\frac{\|{\rot}'f'\|_\infty}{\mu'}=
\frac{\|{\rot}f\|_\infty L^2}{\nu\mu}=D.
\end{equation}
}
\end{remark}

The next lemma is similar to the main estimate for the space analyticity radius
in~\cite{Kuk}.
\begin{lemma}\label{L: analytic sol}
Suppose that $F$ is a restriction to $\Omega$ $($that is, $y=0$$)$
of a bounded $x$-periodic analytic function
$F(x+iy)+iG(x+iy)$ in the region $|y|\le\delta_F$ and
\begin{equation}\label{MF}
M_F^2=\sup_{x\in\Omega, \ |y|\le\delta_F}
(F(x+iy)^2+G(x+iy)^2).
\end{equation}
Let $p\ge3/2$ and let
$$
t_0=\frac{M_{2p}^2}{CM_F^2/\mu}\,.
$$
Here $($and throughout$)$ $C$ is a sufficiently large universal constant
and $M_{2p}\ge\|\omega_0\|_{L_{2p}}$. Then the solution $\omega(t)$
is analytic for $t>0$ and for $0<t\le t_0$ the space analyticity
radius of  $\omega(t)$ is greater than
$$
\delta(t)=
\min\left(
\frac{t^{1/2}}C,
\frac1{Cpt^{(2p-3)/4p}M_{2p}},
\frac1{Cpt^{(2p-3)/(4p+6)}M_{2p}^{2p/(2p+3)}},
\frac1{pt^{1/2}M_{2p}}, \delta_F
\right).
$$
\end{lemma}
\begin{proof}
We solve (\ref{DVE}) by a sequence of approximating solutions (see~\cite{Kato},
\cite{Kuk}).
We set $u^{(0)}=0$ and $\omega^{(0)}=0$. Then for $\omega^{(n)}$,
$u^{(n)}$ we have the equation
\begin{equation}\label{Iter}
\aligned
\partial_t\omega^{(n)}-\Delta\omega^{(n)}+
u^{(n-1)}\cdot\nabla\omega^{(n)}+\mu \omega^{(n)}=F\\
\omega^{(n)}(0)=\omega_0=\rot u_0, \qquad
u^{(n)}=\nabla^\perp\Delta^{-1}\omega^{(n)}.
\endaligned
\end{equation}

 The solutions $\omega^{(n)}$
and $u^{(n)}$ for $t>0$ have analytic extensions
$\omega^{(n)}+i\theta^{(n)}$ and $u^{(n)}+iv^{(n)}$ and since the
system~(\ref{Iter}) is linear, their analyticity radius is at
least $\delta_F$. They satisfy the equation
$$
\partial_t(\omega^{(n)}+i\theta^{(n)})-\Delta(\omega^{(n)}+i\theta^{(n)})+
(u^{(n-1)}+iv^{(n-1)})\cdot\nabla(\omega^{(n)}+i\theta^{(n)})
+\mu (\omega^{(n)}+i\theta^{(n)})=F+iG,
$$
or, equivalently, the system
\begin{equation}\label{Comp-syst}
\aligned
&\partial_t\omega^{(n)}-\Delta\omega^{(n)}+\mu\omega^{(n)}+
u^{(n-1)}\cdot\nabla\omega^{(n)}-v^{(n-1)}\cdot\nabla\theta^{(n)}&=F,\\
&\partial_t\theta^{(n)}-\Delta\theta^{(n)}+\mu\theta^{(n)}+
u^{(n-1)}\cdot\nabla\theta^{(n)}+v^{(n-1)}\cdot\nabla\omega^{(n)}&=G,
\endaligned
\end{equation}
where, as before, $u^{(n)}=\nabla^\perp\Delta^{-1}\omega^{(n)}$,
$v^{(n)}=\nabla^\perp\Delta^{-1}\theta^{(n)}$, and
 the differential operators are taken with respect to $x$.
In view of the analyticity of the solutions we have the Cauchy--Riemann
equations
\begin{equation}\label{Cauchy-Riemann}
\aligned
\frac{\partial\omega^{(n)}}{\partial y_j}\,&=\,-
\frac{\partial\theta^{(n)}}{\partial x_j}\,,\\
\frac{\partial\omega^{(n)}}{\partial x_j}\,&=\,
\frac{\partial\theta^{(n)}}{\partial y_j}\,,\quad j=1,2,
\endaligned
\end{equation}
and the similar equations for $u^{(n)}$ and $v^{(n)}$.

Let $\ve>0$.
We consider the functional
\begin{equation}\label{psin}
\psi_n(t)=\int_0^1\int_{\Omega}
\bigl(\omega^{(n)}(x,\alpha ts,t)^2+\theta^{(n)}(x,\alpha ts,t)^2+
\ve\bigr)^pdxds.
\end{equation}
We also set
$$
Q_n(x,s,t)=\omega^{(n)}(x,\alpha ts,t)^2+\theta^{(n)}(x,\alpha ts,t)^2+\ve.
$$
Here $t\in\mathbb{R}^+$ and $\alpha\in \mathbb{R}^2$. The combination
$\alpha ts$ will play the role of the variable $y$;  $p\ge3/2$, and
$\ve>0$ is arbitrary.

We differentiate $\psi_n(t)$ taking  into account~(\ref{Comp-syst})
and use the Cauchy--Riemann equations~(\ref{Cauchy-Riemann})
 to handle the derivatives with respect to $y$. We obtain
 \begin{equation}\label{I0-I4}
 \frac1{2p}\partial_t\psi_n(t)+I_0=I_1+I_2+I_3+I_4,
\end{equation}
where
$$
\aligned
I_0=\int_0^1\int_\Omega Q_n^{p-1}\bigl(|\nabla\omega^{(n)}|^2+
|\nabla\theta^{(n)}|^2+\mu(\omega^{(n)})^2+\mu(\theta^{(n)})^2\bigr)dxds+\\
+2(p-1)\int_0^1\int_\Omega Q_n^{p-2}\bigl(\omega^{(n)}\nabla\omega^{(n)}+
\theta^{(n)}\nabla\theta^{(n)}\bigr)^2dxds,
\endaligned
$$
and
$$
\aligned
&I_1=\int_0^1\int_\Omega Q_n^{p-1}\bigl(-\omega^{(n)}\nabla\theta^{(n)}+
\theta^{(n)}\nabla\omega^{(n)}\bigr)\cdot\alpha s\,dxds,\\
&I_2=\int_0^1\int_\Omega Q_n^{p-1}\bigl(\omega^{(n)}\nabla\omega^{(n)}+
\theta^{(n)}\nabla\theta^{(n)}\bigr)\cdot u^{(n-1)}dxds,\\
&I_3=\int_0^1\int_\Omega Q_n^{p-1}\bigl(-\omega^{(n)}\nabla\theta^{(n)}+
\theta^{(n)}\nabla\omega^{(n)}\bigr)\cdot v^{(n-1)}dxds,\\
&I_4=\int_0^1\int_\Omega Q_n^{p-1}\bigl(\omega^{(n)}F+
\theta^{(n)}G)dxds.
\endaligned
$$
The arguments of $Q_n$ are $x,s,t$, and the arguments of $\omega^{(n)}$,
$\theta^{(n)}$, $u^{(n)}$, and $v^{(n)}$ are $x$, $\alpha ts$, and $t$.

For an arbitrary $\eta>0$ we have
\begin{equation}\label{I1}
\aligned
I_1&\le\eta
\int_0^1\int_\Omega Q_n^{p-1}\bigl(|\nabla\omega^{(n)}|^2+
|\nabla\theta^{(n)}|^2\bigr)dxds+\\
&C_\eta\int_0^1\int_\Omega Q_n^{p-1}\bigl((\omega^{(n)})^2+
(\theta^{(n)})^2\bigr)|\alpha|^2s^2dxds\le
\eta I_0+C_\eta|\alpha|^2\psi_n(t).
\endaligned
\end{equation}
Next,
\begin{equation}\label{I2}
I_2=\frac1{2p}\int_0^1\int_\Omega \nabla Q_n^{p}\cdot u^{(n-1)}dxds=0.
\end{equation}
For $I_3$ we have
\begin{equation}\label{I3-1}
I_3\le \eta I_0+C_\eta\int_0^1\int_\Omega
Q_n^{p}\,|v^{(n-1)}|^2dxds\le
\eta I_0+C_\eta I_3'I_3'',
\end{equation}
where
$$
I_3'=\left(\int_0^1\int_\Omega Q_n(x,s,t)^{p^2/(p-1)}
dxds\right)^{(p-1)/p}\!,
\,
I_3''=\sum_{j=1}^2\left(\int_0^1\int_\Omega |v_j^{(n-1)}
(x,\alpha t s, t)|^{2p}dxds\right)^{1/p}\!.
$$
We write $I_3'$ as follows
$$
I_3'=\|Q_n^{p/2}\|_{L_\beta(\Omega_0)}^2,
\qquad\Omega_0=\Omega\times[0,1]\subset\mathbb{R}^3,
\qquad\beta=2p/({p-1}),\qquad 2\le\beta\le6,
$$
and use in $\Omega_0$ the Gagliardo--Nirenberg inequality
$$
\|A\|_{L_\beta(\Omega_0)}\le
C\|A\|_{L_2(\Omega_0)}^{3/\beta-1/2}
\|\nabla_{x,s}\,A\|_{L_2(\Omega_0)}^{3/2-3/\beta}+
C\|A\|_{L_2(\Omega_0)}
$$
for $A=A(x,s)=Q_n^{p/2}(x,s,t)$. We have
$$
\aligned
&\|\nabla_{x,s}\,A\|_{L_2(\Omega_0)}^2=
\|\nabla_{x,s}\,Q_n^{p/2}\|_{L_2(\Omega_0)}^2=\\
&p^2\int_0^1\int_\Omega
Q_n^{p-2}\bigl((\omega^{(n)}\nabla\omega^{(n)}
+\theta^{(n)}\nabla\theta^{(n)})^2+ t^2
(\theta^{(n)}\alpha\cdot\nabla\omega^{(n)}-
\omega^{(n)}\alpha\cdot\nabla\theta^{(n)})^2\bigr)dxds\le\\
&\qquad\qquad \le Cp^2(1+|\alpha|^2t^2)I_0.
\endaligned
$$
Hence,
$$
\|\nabla_{x,s}\,A\|_{L_2(\Omega_0)}^{3/2-3/\beta}=
\|\nabla_{x,s}\,Q_n^{p/2}\|_{L_2(\Omega_0)}^{3/2p}\le
C(1+|\alpha|^2t^2)^{3/4p}I_0^{3/4p}.
$$
Next, $\|Q_n^{p/2}\|_{L_2(\Omega_0)}^2=\psi_n(t)$,
$$
\|A\|_{L_2(\Omega_0)}^{3/\beta-1/2}=
\|Q_n^{p/2}\|_{L_2(\Omega_0)}^{(2p-3)/2p}=\psi_n(t)^{(2p-3)/4p}
$$
and
\begin{equation}\label{I3'}
I_3'\le C(1+|\alpha|^2t^2)^{3/2p}I_0^{3/2p}\psi_n(t)^{(2p-3)/2p}+C\psi_n(t).
\end{equation}
We now consider $I_3''$. Since $v_j^{(n-1)}(x,0,t)=0$
(the solution restricted to $y=0$ is real-valued), we have
(using the Cauchy--Riemann equations for $v_j$)
$$
\aligned
|v_j^{(n-1)}(x,\alpha ts,t)|=&\left|\sum_{k=1}^2\alpha_kts
\int_0^1\partial_{y_k} v_j^{(n-1)}(x,\alpha ts\tau,t)d\tau\right|=\\
&\left|\sum_{k=1}^2\alpha_kts
\int_0^1\partial_{k} u_j^{(n-1)}(x,\alpha ts\tau,t)d\tau\right|.
\endaligned
$$
Then
$$
\aligned
I_3''=\sum_{j=1}^2\left(\int_0^1\int_\Omega\left|\sum_{k=1}^2
\alpha_kts\int_0^1\partial_ku_j^{(n-1)}(x,\alpha ts\tau,t)d\tau\right|^{2p}
dxds\right)^{1/p}\le\\
C|\alpha|^2t^2\left(\int_0^1\int_\Omega\int_0^1\left|
\nabla u^{(n-1)}(x,\alpha ts\tau,t)\right|^{2p}s^{2p}d\tau
dxds\right)^{1/p}=\\
C|\alpha|^2t^2\left(\int_0^1s^{2p}ds\int_0^1d\tau\int_\Omega\left|
\nabla u^{(n-1)}(x,\alpha ts\tau,t)\right|^{2p}
dx\right)^{1/p}.
\endaligned
$$
Since $u=\nabla^\perp\Delta^{-1}\omega$, we have (see~\cite{G-T}, \cite{Yud})
$$
\left(\int_\Omega|\nabla u(x)|^{2p}dx\right)^{1/2p}=
\|\nabla\nabla^\perp\Delta^{-1}\omega\|_{L_{2p}}\le
\|\Delta^{-1}\omega\|_{W_{2p}^2}\le
Cp\|\omega\|_{L_{2p}}\,.
$$
Therefore
$$
\aligned
I_3''\le
Cp^2|\alpha|^2t^2\left(\int_0^1\int_0^1\int_\Omega\left|
\omega^{(n-1)}(x,\alpha ts\tau,t)\right|^{2p}
dx\, s^{2p}dsd\tau\right)^{1/p}\le\\
Cp^2|\alpha|^2t^2\left(\int_0^1\int_\Omega Q_{n-1}^pdxds\right)^{1/p}\le
Cp^2|\alpha|^2t^2\psi_{n-1}(t)^{1/p},
\endaligned
$$
where we have used
$\int_0^1\int_0^1h(s\tau)s^{2p}dsd\tau\le(2p)^{-1}\int_0^1h(s)ds$.
Combining this with~(\ref{I3-1}) and (\ref{I3'}) we obtain
\begin{equation}\label{I3-final}
\aligned
&I_3\le\\ &\eta'I_0+C_{\eta'}p^2|\alpha|^2t^2(1+|\alpha|^2t^2)^{3/2p}
I_0^{3/2p}\psi_n(t)^{(2p-3)/2p}\psi_{n-1}(t)^{1/p}+
Cp^2|\alpha|^2t^2\psi_{n-1}(t)^{1/p}
\psi_n(t)\le\\
&\eta I_0+C_{\eta}p^2(|\alpha|t)^{4p/(2p-3)}(1+|\alpha|^2t^2)^{3/(2p-3)}
\psi_{n-1}(t)^{2/(2p-3)}\psi_{n}(t)+Cp^2|\alpha|^2t^2\psi_{n-1}(t)^{1/p}
\psi_n(t).
\endaligned
\end{equation}
Finally, we estimate $I_4$:
\begin{equation}\label{I4}
\aligned
I_4\le\int_0^1\int_\Omega Q_n^{p-1}\bigl((\omega^{(n)})^2
+\theta^{(n)})^2)\eta\mu+(F^2+G^2)/(4\eta\mu)\bigr)dxds\le\\
\eta I_0+C_\eta(M_F^2/\mu)\int_0^1\int_\Omega Q_n^{p-1}dxds
\le\eta I_0+C_\eta(M_F^2/\mu)\psi_n(t)^{(p-1)/p},
\endaligned
\end{equation}
where $M_F$ is defined in~(\ref{MF}).

Taking $\eta>0$ sufficiently small we infer from~(\ref{I0-I4}),
(\ref{I1}), (\ref{I2}), (\ref{I3-final}), (\ref{I4})
$$
\aligned
\partial_t\psi_n(t)\le Cp|\alpha|^2\psi_n(t)+Cp^3|\alpha|^{4p/(2p-3)}
t^{4p/(2p-3)}\psi_{n-1}(t)^{2/(2p-3)}\psi_n(t)+\\
Cp^3|\alpha|^{(4p+6)/(2p-3)}
t^{(4p+6)/(2p-3)}\psi_{n-1}(t)^{2/(2p-3)}\psi_n(t)+\\
Cp^3|\alpha|^2t^2\psi_{n-1}(t)^{1/p}\psi_n(t)+
Cp\psi_n(t)^{(p-1)/p}M_F^2/\mu,
\endaligned
$$
where
$$
\psi_n(0)=\int_\Omega(\omega_0(x)^2+\ve)^pdx.
$$
We set $\vphi_n(t)=\psi_n(t)^{1/2p}$ and obtain the differential inequality for
$\vphi_n$:
$$
\aligned
\partial_t\vphi_n(t)\le C|\alpha|^2\vphi_n(t)+Cp^2|\alpha|^{4p/(2p-3)}
t^{4p/(2p-3)}\vphi_{n-1}(t)^{4p/(2p-3)}\vphi_n(t)+\\
Cp^2|\alpha|^{(4p+6)/(2p-3)}
t^{(4p+6)/(2p-3)}\vphi_{n-1}(t)^{4p/(2p-3)}\vphi_n(t)+\\
Cp^2|\alpha|^2t^2\vphi_{n-1}(t)^{2}\vphi_n(t)+
C\vphi_n(t)^{-1}M_F^2/\mu.
\endaligned
$$

We now use the Gronwall-type Lemma~\ref{L:Kuk-Gronwall} from \cite{Kuk} below
and see that $\vphi(t)\le2\vphi(0)$ on the time interval specified in~(\ref{tmin}),
(\ref{A1-A5}),
and letting $\ve\to0$ we obtain
$$
\int_0^1\int_\Omega\left(\omega^{(n)}(x,\alpha ts,t)^2+\theta^{(n)}(x,\alpha ts,t)^2
\right)^pdxds\le2^{2p}M_{2p}^{2p},
$$
for $t\ge0$, $|\alpha|t\le\delta_F$
and
\begin{equation}\label{tmin}
t\le\min(A_1,A_2,A_3,A_4,A_5),
\end{equation}
where
\begin{equation}\label{A1-A5}
\aligned
&A_1=\frac1{C|\alpha|^2},\\
&A_2=\frac1{Cp^{(4p-6)/(6p-3)}|\alpha|^{4p/(6p-3)}M_{2p}^{4p/(6p-3)}},\\
&A_3=\frac1{Cp^{(4p-6)/(6p+3)}|\alpha|^{(4p+6)/(6p+3)}M_{2p}^{4p/(6p+3)}},\\
&A_4=\frac1{Cp^{2/3}|\alpha|^{2/3}M_{2p}^{2/3}},\\
&A_5=\frac{M_{2p}^2}{CM_F^2\mu^{-1}}\,.
\endaligned
\end{equation}

We now set
\begin{equation}\label{t0}
t_0=\frac{M_{2p}^2}{CM_F^2/\mu}\,.
\end{equation}
Then the condition
$$
t\le\min(A_1,A_2,A_3,A_4)
$$
can be written in terms of $y=\alpha t$ as follows
\begin{equation}\label{|y|}
|y|\le\min\left(
\frac{t^{1/2}}C,
\frac1{Cpt^{(2p-3)/4p}M_{2p}},
\frac1{Cpt^{(2p-3)/(4p+6)}M_{2p}^{2p/(2p+3)}},
\frac1{pt^{1/2}M_{2p}}
\right).
\end{equation}
Now for $t_0$ defined in~(\ref{t0}) and
\begin{equation}\label{delta(t)}
\delta(t)=
\min\left(
\frac{t^{1/2}}C,
\frac1{Cpt^{(2p-3)/4p}M_{2p}},
\frac1{Cpt^{(2p-3)/(4p+6)}M_{2p}^{2p/(2p+3)}},
\frac1{pt^{1/2}M_{2p}}, \delta_F
\right)
\end{equation}
we have for $0<t\le t_0$ and $|y|\le\delta(t)$
$$
\int_0^1\int_\Omega\bigl(\omega^{(n)}(x,sy,t)^2
+\theta^{(n)}(x,sy,t)^2\bigr)dxds\le2^{2p}M_{2p}^{2p}
$$
for all integer $n\ge1$. Therefore for any $y\in\mathbb{R}^2$ with $|y|=1$
this gives that
$$
\int_0^{\delta(t)}\int_\Omega\bigl(\omega^{(n)}(x,sy,t)^2
+\theta^{(n)}(x,sy,t)^2\bigr)dxds\le2^{2p}\delta(t)M_{2p}^{2p}
$$
and since $\int_0^\delta f(sy)ds\le B$, $|y|=1$ implies
$\int_{|y|\le\delta}f(y)dy\le2\pi\delta B$,
we obtain
$$
\int_{|y|\le\delta(t)}\int_\Omega\bigl(\omega^{(n)}(x,y,t)^2
+\theta^{(n)}(x,y,t)^2\bigr)dxds\le2\pi2^{2p}\delta(t)^2M_{2p}^{2p}.
$$
This estimate is uniform in $n$ and as in~\cite{Gr-K}, \cite{Kuk} we obtain the
existence of an analytic solution of~(\ref{DVE}) with analyticity radius
satisfying~(\ref{delta(t)}). The proof is complete.
\end{proof}

\begin{lemma}\label{L:Kuk-Gronwall} {\rm(See~\cite{Kuk}.)}
Let  $y_n(t)\in C^1[0,T]$ be a sequence of non-negative functions satisfying
$y_0(t)\le M$ for $0\le t\le T$, and $y_n(0)\le M$ for $n\ge1$.
Suppose that on the interval $0\le t\le T$
$$
 \partial_ty_n(t)\le\sum_{j=1}^NK_jt^{\alpha_j}y_n(t)^{\beta_j}y_{n-1}(t)^{\gamma_j},
 $$
 where $K_j>0$, $\alpha_j>-1$,  $\beta_j\in\mathbb{R}$, and  $\gamma_j\ge0$
 are given constants. Then $y_n(t)\le 2M$ for all $n=0,1,2,\dots$ provided that
 $$
 0\le t\le\min\left(T, \min_{j=1,\dots,N}\left(\frac{\alpha_j+1}
 {NK_j2^{\beta_j^++\gamma_j}M^{\beta_j+\gamma_j-1}}
 \right)^{1/(\alpha_j+1)}\right),
 $$
 where $\beta^+=\max(\beta,0)$.
 \end{lemma}

 We can now state the main result of this section.
\begin{theorem}\label{T:Analyticity _in_L_infty}
The solutions on the 2D space-periodic damped-driven Navier--Stokes system~$(\ref{DNS})$
lying on the global attactor $\mathcal{A}$ are analytic with
space analyticity radius $l_a$ satisfying the lower bound
\begin{equation}\label{est_for_la}
l_a\ge\frac{|\Omega|^{1/2}}{CD^{1/2}(1+\log D)^{1/2}},\qquad
\text{where}\qquad
D=\frac{\|\rot f\|_\infty|\Omega|}{\mu\nu}\,.
\end{equation}
\end{theorem}
\begin{proof}
We first observe that~(\ref{est_for_la}) is equivalent
to the estimate
\begin{equation}\label{est_for_la_dimless}
l_a\ge\frac{1}{CD^{1/2}(1+\log D)^{1/2}}
\end{equation}
for the equation written in dimensionless form.

Next, by Young's inequality
$$
\aligned
&pt^{(2p-3)/4p}M\le CpM^{4p/(4p-3)}t^{1/2}+t^{-1/2},\\
&pt^{(2p-3)/(4p+6)}M^{2p/(2p+3)}\le CpM^{}t^{1/2}+t^{-1/2}.
\endaligned
$$
Hence,
the estimate~(\ref{delta(t)}) can be written as follows
\begin{equation}\label{delta(t)3}
\delta(t)\ge
\min\left(
\frac{t^{1/2}}C,
\frac1{Cpt^{1/2}\bigl(M_{2p}^{4p/(4p-3)}+M_{2p}\bigr)}, \delta_F
\right).
\end{equation}
The solutions lying on the attractor are bounded in $L_{2p}$:
$$
\|\omega(t)\|_{L_{2p}}\le M_{2p},
\qquad
\M_{2p}\le CM_\infty.
$$
Setting
$$
p=C(1+\log M_\infty)
$$
we see that
$$
\left.p\bigl(M_{2p}^{4p/(4p-3)}+M_{2p}\bigr)\right\vert_{p=C(1+\log M_\infty)}
\le C(1+\log M_\infty)M_\infty
$$
and therefore
$$
\delta(t)\ge
\min\left(
\frac{t^{1/2}}C\,,\,
\frac1{C(1+\log M_\infty)M_\infty t^{1/2}}\,, \, \delta_F
\right).
$$
At the moment of time
$$
t^*=\frac1{C(1+\log M_\infty)M_\infty}\,,
$$
which for sufficiently large $M_\infty$ (the case of our interest) is smaller
than $t_0$ defined in~(\ref{t0}) (the details are given below) we have
$$
\delta(t^*)\ge
\frac1{CM_\infty^{1/2}(1+\log M_\infty)^{1/2}}\,.
$$
Since $M_\infty\le D$ (see (\ref{Dimless-est})), it follows that
$$
\delta(t^*)\ge
\frac1{CD^{1/2}(1+\log D)^{1/2}}\,.
$$
By the invariance property of the attractor we see that on the attractor  the above estimate holds for all $t^*$, which proves~(\ref{est_for_la_dimless}).

To complete the proof it remains to show that
\begin{equation}\label{t0t*}
\frac1{C(1+\log D)D}=t^*\le t_0=\frac{D^2}{CM_{F'}^2/\mu'}\,,
\end{equation}
where in the expression for $t_0$ we reverted to the prime notation
for the dimensionless damping coefficient $\mu'$ and the forcing $F'$.
We relate the forcing term and its analytic extension by the equality
$$
M_{F'}=K\|{\rot}'f'\|_\infty,\qquad K=K(F,\delta_F).
$$
Recalling that $f'=(L^3/\nu^2)f$, $\mu'=(L^2/\nu)\mu$, and $x'=(1/L)x$
we see that
$$
t_0=\frac{\nu}{CK^2\mu L^2}\,.
$$
Hence, (\ref{t0t*}) goes over to the condition
$$
C(1+\log D)\ge\frac{K^2\mu^2}{\|\rot f\|_\infty}\,,
$$
which is obviously satisfied for all sufficiently small $\nu>0$.
The proof is complete.
\end{proof}

\setcounter{equation}{0}
\section{Concluding remarks}\label{S:Conclusion}

We have shown that the solutions lying on the attractor of the 2D space-periodic
damped-driven Navier--Stokes system,
the Stommel--Charney barotropic model of ocean circulation without rotation,
 with analytic forcing have space
ana\-lyti\-city radius which up to a logarithmic term coincides with the small
scale estimates both  in terms of the sharp bounds for the fractal dimension of
 the global attractor, and in terms of the spatial lattice of determining nodes.
 The derivation of this lower bound for the analyticity radius essentially uses the techniques developed in~\cite{Kuk}.

\section*{Acknowledgments}

A.A.I. would like to thank the warm hospitality of the Mathematics Department
at the University of California, Irvine, where this work was done.

 This work was supported in part by the US Civilian Research and
Development Foundation, grant no.~RUM1-2654-MO-05 ( A.A.I. and E.S.T.), by the
Russian Foundation for Fundamental Research, grants~nos.~06-01-00096
and 05-01-00429,  and by the
RAS Programme no.1 `Modern problems of theoretical mathematics' (A.A.I.). The work of E.S.T. was supported in part by the National
Science Foundation, grant no.~DMS-0504619, the ISF grant no. 120/6, and the
BSF grant no. 2004271.

\bibliographystyle{amsplain}

\end{document}